\overfullrule=0pt
\centerline {\bf A class of equations with three solutions}\par
\bigskip
\bigskip
\centerline {BIAGIO RICCERI}\par
\bigskip
\bigskip
{\bf Abstract:} Here is one of the results obtained in this paper: Let $\Omega\subset {\bf R}^n$ be a smooth bounded domain, let $q>1$, with $q<{{n+2}\over {n-2}}$ if
$n\geq 3$ and let $\lambda_1$ be the first eigenvalue of the problem
$$\cases{-\Delta u=\lambda u & in $\Omega$ \cr & \cr u=0 & on $\partial\Omega$\ .\cr}$$
Then, for every $\lambda>\lambda_1$ and for every convex set $S\subseteq L^{\infty}(\Omega)$ dense in $L^2(\Omega)$, there exists
$\alpha\in S$ such that the problem
$$\cases{-\Delta u=\lambda(u^+-(u^+)^q)+\alpha(x) & in $\Omega$ \cr & \cr u=0 & on $\partial\Omega$\cr}$$
has at least three weak solutions, two of which are global minima in $H^1_0(\Omega)$ of the functional
$$u\to {{1}\over {2}}\int_{\Omega}|\nabla u(x)|^2dx-\lambda\int_{\Omega}\left ({{1}\over {2}}|u^+(x)|^2-{{1}\over {q+1}}|u^+(x)|^{q+1}\right )dx-\int_{\Omega}\alpha(x)u(x)dx\ $$
where $u^+=\max\{u,0\}$.
\bigskip
\bigskip
\bigskip
\bigskip
{\bf 1. Introduction}\par
\bigskip
There is no doubt that the study of nonlinear PDEs lies in the core of Nonlinear Analysis. In turn, one of the most studied topics concerning nonlinear PDEs is
the multiplicity of solutions. On the other hand, the study of the global minima of integral functionals is essentially the central subject of the Calculus of Variations.
In the light of these facts, it is hardly understable why the number of the known results on multiple global minima of integral functionals is extremely low. Certainly, this is not due
to a lack of intrinsic mathematical interest. Probably, the reason could reside in the fact that there is not an abstract tool which has the same popularity as the one that, for instance,
the Lyusternik-Schnirelmann theory and the Morse theory have in dealing with multiple solutions for nonlinear PDEs.\par
\smallskip
 Abstract results on the multiplicity of
global minima, however, are already present in the literature. We mainly allude to the result first obtained in [1] and then extended in [2] and [5] which ensures the existence of at least
two global minima provided that a strict minimax inequality holds. We already have obtained a variety of applications upon different ways of checking the required strict inequality 
([3], [4], [6]).\par
\smallskip
The aim of the present paper is to establish an application of Theorem 1 of [7] which is itself an application of the main result in [5]. Precisely, we first establish a general result which ensures the existence of three solutions for a certain equation provided that another related one has no non-zero solutions (Theorem 1). Then, we present an application to nonlinear elliptic equations (Theorem 2).\par
\bigskip
{\bf 2. Results}\par
\bigskip
In the sequel, $(X,\|\cdot\|_X)$ is a reflexive real Banach space, $(Y,\langle\cdot,\cdot\rangle_Y)$ is a real Hilbert space, $I, \psi:X\to {\bf R}$ are two $C^1$ functionals, with $I(0)=\psi(0)=0$
and $\sup_{\bf R}\psi>0$, 
$\varphi:X\to Y$ is a $C^1$ operator, with $\varphi(0)=0$. 
For each fixed $y\in Y$, we denote by $\partial_x\langle \varphi(\cdot),y\rangle_Y$ the derivative
of the functional $x\to \langle \varphi(x),y\rangle_Y$. Clearly, one has
$$\partial_x\langle\varphi(x),y\rangle_Y (u)=\langle \varphi'(x)(u),y\rangle_Y$$
for all $x, u\in X$. \par
\smallskip
We say that $I$ is coercive if $\lim_{\|x\|_X\to +\infty}I(x)=+\infty$. We also say that $I'$ admits a continuous inverse on $X^*$ if there exists
a continuous operator $T:X^*\to X$ such that $T(I'(x))=x$ for all $x\in X$.\par
\smallskip
Here is our abstract result:\par
\medskip
THEOREM 1. {\it Let $I$ be weakly lower semicontinuous and coercive, and let $I'$ admit a continuous inverse on $X^*$. Moreover,
assume that the operators $\varphi'$ and $\psi'$ are compact and that
$$\lim_{\|x\|_X\to +\infty}{{\langle\varphi(x),y\rangle_Y}\over {I(x)}}=0\eqno{(1)}$$
for all $y$ in a convex and dense set $T\subseteq Y$.
Set
$$\theta^*:=\inf_{x\in \psi^{-1}(]0,+\infty[)}{{I(x)}\over {\psi(x)}}\ ,$$
$$\tilde\theta:=\cases{\liminf_{x\in \psi^{-1}(]0,+\infty[), \|x\|_X\to +\infty}{{I(x)}\over {\psi(x)}} & if $\psi^{-1}(]0,+\infty[)$\hskip 3pt
is\hskip 3pt unbounded\cr & \cr +\infty & otherwise\cr}$$
and assume that
$$\theta^*<\tilde\theta\ .$$
Then, for each  $\lambda\in ]\theta^*,\tilde\theta[$, with $\lambda\geq 0$, either
the equation
$$I'(x)=-\partial_x\langle \varphi(x),\varphi(x)\rangle_Y+\lambda\psi'(x)$$
has a non-zero solution, or,
 for each convex set $S\subseteq T$ dense in $Y$, there exists $\tilde y\in S$ such that
the equation
$$I'(x)=\partial_x\langle \varphi(x),\tilde y\rangle_Y+\lambda\psi'(x)$$
has at least three solutions, two of which are global minima in $X$ of the functional $$x\to I(x)-\langle\varphi(x),\tilde y\rangle_Y-
\lambda\psi(x)\ .$$}\par
\medskip
As it was said in the Introduction, the main tool to prove Theorem 1 is a result recently obtained in [7].
For reader's convenience, we now recall its statement:\par
\medskip
THEOREM A ([7], Theorem 1). - {\it Let $X, E$ be two real reflexive Banach spaces
 and let $\Phi:X\times E\to {\bf R}$ be a $C^1$ functional satisfying the following
conditions:\par
\noindent
$(a)$\hskip 5pt 
the functional $\Phi(x,\cdot)$ is quasi-concave for all $x\in X$ and the functional $-\Phi(x_0,\cdot)$ is coercive for some $x_0\in X$;\par
\noindent
$(b)$\hskip 5pt  there exists a convex set $T\subseteq E$ dense in $E$, such that, for each $y\in T$, the functional
$\Phi(\cdot,y)$ is weakly lower semicontinuous, coercive and satisfies the Palais-Smale condition\ .\par
Then, either the system
$$\cases {\Phi'_x(x,y)=0 \cr & \cr \Phi'_y(x,y)=0\cr}$$
has a solution $(x^*,y^*)$ such that
$$\Phi(x^*,y^*)=\inf_{x\in X}\Phi(x,y^*)=\sup_{y\in E}\Phi(x^*,y)\ ,$$
or, for every convex set $S\subseteq T$ dense in $E$, there exists $\tilde y\in S$ such that
equation
$$\Phi'_x(x,\tilde y)=0$$ has at least three solutions, two of which are global minima in $X$ of the functional $\Phi(\cdot,\tilde y)$.}\par
\medskip
{\it Proof of Theorem 1}. Fix $\lambda\in ]\theta^*,\tilde\theta[$, with $\lambda\geq 0$. Assume that
the equation
$$I'(x)=-\partial_x\langle \varphi(x),\varphi(x)\rangle_Y+\lambda\psi'(x)$$
has no non-zero solution. Fix a convex set $S\subseteq T$ dense in $Y$. We have to show that there exists $\tilde y\in S$ such that the equation
$$I'(x)=\partial_x\langle \varphi(x),\tilde y\rangle_Y+\lambda\psi'(x)$$
has at least three solutions, two of which are global minima in $X$ of the functional $x\to I(x)-\langle\varphi(x),\tilde y\rangle_Y-
\lambda\psi(x)$. To this end, let us apply Theorem A.
Consider the functional
$\Phi:X\times Y\to {\bf R}$ defined by
$$\Phi(x,y)=I(x)-{{1}\over {2}}\|y\|_Y^2-\langle\varphi(x),y\rangle-\lambda\psi(x)$$
for all $(x,y)\in X\times Y$. 
Of course, $\Phi$ is $C^1$ and, for each $x\in X$, $\Phi(x,\cdot)$ is concave and $-\Phi(x,\cdot)$ is coercive.
Fix $y\in T$. Let us show that the operator $\partial_x\langle\varphi(\cdot),y\rangle$ is compact. To this end, let $\{x_n\}$ be a bounded
sequence in $X$. Since $\varphi'$ is compact, up to a subsequence, $\{\varphi'(x_n)\}$ converges in ${\cal L}(X,Y)$ to
some $\eta$.  That is
$$\lim_{n\to \infty}\sup_{\|u\|_X=1}\|\varphi'(x_n)(u)-\eta(u)\|_Y=0\ .$$
On the other hand, we have
$$\sup_{\|u\|_X=1}|\partial_x\langle\varphi(x_n),y\rangle_Y (u)-\langle\eta(u),y\rangle_Y|=\sup_{\|u\|_X=1}|\langle\varphi'(x_n)(u),y\rangle_Y-
\langle\eta(u),y\rangle_Y|$$
$$\leq \sup_{\|u\|_X=1}\|\varphi'(x_n)(u)-\eta(u)\|_Y\|y\|_Y$$
and so the sequence $\{\partial_x\langle\varphi(x_n),y\rangle_Y (\cdot)\}$ converges in $X^*$ to $\eta(\cdot)(y)$. Then, since $\psi'$ is
compact, the operator $\partial_x\langle\varphi(\cdot),y\rangle_Y+\lambda\psi'(\cdot)$ is compact too. From this, it
follows that $\langle\varphi(\cdot),y\rangle_Y+\lambda\psi(\cdot)$ is sequentially weakly continuous ([8], Corollary 41.9). If $\|x\|_X$ is large enough,
we have $I(x)>0$ and so we can write
$$\Phi(x,y)=I(x)\left (1-{{{{1}\over {2}}\|y\|_Y^2+\langle\varphi(x),y\rangle_Y+\lambda\psi(x)}\over {I(x)}}\right )\ .\eqno{(2)}$$
In view of $(1)$, we also have
$$\liminf_{\|x\|_X\to +\infty}\left (1-{{{{1}\over {2}}\|y\|_Y^2+\langle\varphi(x),y\rangle_Y+\lambda\psi(x)}\over {I(x)}}\right )=
1- \limsup_{\|x\|_X\to +\infty}{{\lambda\psi(x)}\over {I(x)}}\ .\eqno{(3)}$$
We claim that
$$\limsup_{\|x\|\to +\infty}{{\lambda\psi(x)}\over {I(x)}}<1\ .\eqno{(4)}$$
This is clear if either $\lambda=0$ or $\limsup_{\|x\|_X\to +\infty}{{\psi(x)}\over {I(x)}}\leq 0$. If $\lambda>0$ and $\limsup_{\|x\|_X\to +\infty}{{\psi(x)}\over {I(x)}}>0$, then $(4)$
is equivalent to
 $$\limsup_{\|x\|_X\to +\infty}{{\psi(x)}\over {I(x)}}<+\infty$$ 
and
$$\lambda<{{1}\over {\limsup_{\|x\|_X\to +\infty}{{\psi(x)}\over {I(x)}}}}\ .\eqno{(5)}$$
But
$${{1}\over {\limsup_{\|x\|_X\to +\infty}}{{\psi(x)}\over {I(x)}}}=\liminf_{x\in \psi^{-1}(]0,+\infty[), \|x\|_X\to +\infty}{{I(x)}\over {\psi(x)}}\ ,$$
and so $(5)$ is satisfied just since $\lambda<\tilde\theta$.
Since $I$ is coercive
and weakly lower semicontinuous, the functional $\Phi(\cdot,y)$ turns out to be coercive, in view of $(2)$, $(3)$, $(4)$, and weakly lower semicontinuous,
in view of the Eberlein-Smulyan theorem. Finally, since $I'$ admits a continuous inverse on
$X^*$, $\Phi(\cdot,y)$ satisfies the Palais-Smale condition in view of Example 38.25 of [8]. Hence, $\Phi$
satisfies the assumptions of Theorem A.  Now,  we claim that there is no solution $(x^*,y^*)$ of the system
$$\cases {\Phi'_x(x,y)=0 \cr & \cr \Phi'_y(x,y)=0\cr}$$
such that 
$$\Phi(x^*,y^*)=\inf_{x\in X}\Phi(x,y^*)\ .$$
Arguing by contradiction, assume that such a $(x^*,y^*)$ does exist.
This amounts to say that 
$$\cases {I'(x^*)=\partial_x\langle\varphi(x^*),y^*\rangle_Y+\lambda\psi'(x^*) \cr & \cr y^*=-\varphi(x^*)\cr}$$
and
$$I(x^*)-\langle\varphi(x^*),y^*\rangle_Y-\lambda\psi(x^*)=\inf_{x\in X}(I(x)-\langle\varphi(x),y^*\rangle_Y-\lambda\psi(x))\ .\eqno{(6)}$$
Therefore
$$I'(x^*)=-\partial_x\langle\varphi(x^*),\varphi(x^*)\rangle_Y+\lambda\psi'(x^*)\ .$$
So, by the initial assumption, we have $x^*=0$ and hence $y^*=0$ (recall that $\varphi(0)=0$). As a consequence, since $I(0)=\psi(0)=0$, $(6)$ becomes
$$\inf_{x\in X}(I(x)-\lambda\psi(x))=0\ .\eqno{(7)}$$
Now, notice that $(7)$ contradicts the fact that $\lambda>\theta^*$. Hence, {\it a fortiori}, the system
$$\cases {\Phi'_x(x,y)=0 \cr & \cr \Phi'_y(x,y)=0\cr}$$
has no solution $(x^*,y^*)$ such that 
$$\Phi(x^*,y^*)=\inf_{x\in X}\Phi(x,y^*)=\sup_{y\in Y}\Phi(x^*,y)$$
and then the existence of $\tilde y\in S$ is directly ensured by Theorem A. \hfill $\bigtriangleup$\par
\medskip
We now present an application of Theorem 1 to a class of nonlinear elliptic equations.\par
\smallskip
Let $\Omega\subset {\bf R}^n$ be a smooth bounded domain. We denote by ${\cal A}$ the class of all Carath\'eodory's functions
$f:\Omega\times {\bf R}\to {\bf R}$ such that, for each $u, v\in H^1_0(\Omega)$, the function $x\to f(x,u(x))v(x)$ lies in $L^1(\Omega)$.
\smallskip
For $f\in {\cal A}$, we consider the Dirichlet problem
$$\cases{-\Delta u=f(x,u) & in $\Omega$ \cr & \cr u=0 & on $\partial\Omega$\ .\cr}$$
As usual,  a weak solution of the problem is any $u\in H^1_0(\Omega)$ such that
$$\int_{\Omega}\nabla u(x)\nabla v(x)dx=\int_{\Omega}f(x,u(x))v(x)dx$$
for all $v\in H^1_0(\Omega)$.\par
\smallskip
For any continuous function $f:{\bf R}\to {\bf R}$, we set $F(\xi)=\int_0^{\xi}f(t)dt$ for all $\xi\in {\bf R}$.\par
\medskip
THEOREM 2. - {\it Let $f, g:{\bf R}\to {\bf R}$ be two continuous functions satisfying the following growth conditions:\par
\noindent
$(a)$\hskip 5pt if $n\leq 3$, one has
$$\lim_{|\xi|\to +\infty}{{|F(\xi)|}\over {\xi^2}}=0\ ;$$
$(b)$\hskip 5pt if $n\geq 2$, there exist $p, q>0$, with $p<{{2}\over {n-2}}$, $q<{{n+2}\over {n-2}}$ if $n\geq 3$, such that
$$\sup_{\xi\in {\bf R}}{{|f(\xi)|}\over {1+|\xi|^p}}<+\infty\ ,$$
$$\sup_{\xi\in {\bf R}}{{|g(\xi)|}\over {1+|\xi|^q}}<+\infty\ .$$
Set
$$\rho:=\limsup_{|\xi|\to +\infty}{{G(\xi)}\over {\xi^2}}\ ,$$
$$\sigma:=\max\left\{\liminf_{\xi\to 0^+}{{G(\xi)}\over {\xi^2}},\liminf_{\xi\to 0^-}{{G(\xi)}\over {\xi^2}}\right\}$$
and assume that
$$\max\{\rho,0\}<\sigma\ .$$
Then, for every $\lambda\in \left ]{{\lambda_1}\over {2\sigma}},{{\lambda_1}\over {2\max\{\rho,0\}}}\right [$ (with the conventions
${{\lambda_1}\over {+\infty}}=0$, ${{\lambda_1}\over {0}}=+\infty$),
 either the problem
$$\cases{-\Delta u=-F(u)f(u)+\lambda g(u) & in $\Omega$ \cr & \cr u=0 & on $\partial\Omega$\cr} \eqno{(8)}$$
has a non-zero weak solution, or,
for every convex set $S\subseteq L^{\infty}(\Omega)$ dense in $L^2(\Omega)$, there exists
$\alpha\in S$ such that the problem
$$\cases{-\Delta u=\alpha(x)f(u)+\lambda g(u) & in $\Omega$ \cr & \cr u=0 & on $\partial\Omega$\cr} \eqno{(9)}$$
has at least three weak solutions, two of which are global minima in $H^1_0(\Omega)$ of the functional
$$u\to {{1}\over {2}}\int_{\Omega}|\nabla u(x)|^2dx-\int_{\Omega}\alpha(x)F(u(x))dx-\lambda\int_{\Omega}G(u(x))dx\ .$$}
\smallskip
PROOF. We are going to apply Theorem 1 taking $X=H^1_0(\Omega)$, $Y=L^2(\Omega)$, with their usual scalar products (that is,
$\langle u,v\rangle_X=\int_{\Omega}\nabla u(x)\nabla v(x)dx$ and $\langle u,v\rangle_Y=\int_{\Omega}u(x)v(x)dx$), $T=L^{\infty}(\Omega)$ and
$$I(u)={{1}\over {2}}\|u\|_X^2\ ,$$
$$\varphi(u)=F\circ u\ ,$$
$$\psi(u)=\int_{\Omega}G(u(x))dx$$
for all $u\in X$. In view of $(b)$, thanks to the Sobolev embedding theorem, the operator $\varphi$ and the functional $\psi$ are
$C^1$, with compact derivative. Moreover, the solutions of the equation
$$I'(u)=-\partial_u\langle\varphi(u),\varphi(u)\rangle_{Y}+\lambda\psi'(u)$$
are weak solutions of $(8)$ and, for each $\alpha\in Y$, the solutions of the equation
$$I'(u)=\partial_u\langle\varphi(u),\alpha\rangle_Y+\lambda\psi'(u)$$
are weak solutions of $(9)$.
Moreover, condition $(1)$ follows readily from $(a)$ which is automatically satisfied when $n\geq 4$ since
$p<{{2}\over {n-2}}$. Now, denote by $\lambda_1$ the first eigenvalue of the Dirichlet problem
$$\cases{-\Delta u=\lambda u & in $\Omega$ \cr & \cr u=0 & on $\partial\Omega$\ .\cr}$$
We claim that
$$\limsup_{\|u\|_X\to +\infty}{{\psi(u)}\over {\|u\|_X^2}}\leq {{\rho}\over {\lambda_1}}\ .\eqno{(10)}$$
Indeed, fix $\nu>\rho$. Then, there exists $\delta>0$ such that
$$G(\xi)\leq \nu\xi^2 \eqno{(11)}$$
for all $\xi\in {\bf R}\setminus [-\delta,\delta]$. Fix $u\in X\setminus \{0\}$. From $(11)$ we clearly obtain
$$\psi(u)\leq \nu\|u\|_Y^2+\hbox {\rm meas}(\Omega)\sup_{[-\delta,\delta]}G\leq
\nu{{\|u\|_X^2}\over {\lambda_1}}+\hbox {\rm meas}(\Omega)\sup_{[-\delta,\delta]}G$$
and so
$$\limsup_{\|u\|_X\to +\infty}{{\psi(u)}\over {\|u\|_X^2}}\leq {{\nu}\over {\lambda_1}}\ .\eqno{(12)}$$
Now, we get $(10)$ passing in $(12)$ to the limit for $\nu$ tending to $\rho$. We also claim that
$$ {{\sigma}\over {\lambda_1}}\leq \sup_{u\in X\setminus \{0\}}{{\psi(u)}\over {\|u\|_X^2}}\ .\eqno{(13)}$$
Indeed, fix $\eta<\sigma$. For instance, let $\sigma=\liminf_{\xi\to 0^+}{{G(\xi)}\over {\xi^2}}$. Then, there exists
$\eta>0$ such that
$$G(\xi)\geq \eta\xi^2 \eqno{(14)}$$
for all $\xi\in [0,\eta]$. Fix any $v\in H^1_0(\Omega)$ such that $\|v\|_X^2=\lambda_1\|v\|_Y^2$ and
$v(\Omega)\subseteq [0,\eta]$. From $(14)$ we obtain
$$\psi(v)\geq\eta\|v\|_Y^2$$
and so
$$\sup_ {u\in X\setminus \{0\}}{{\psi(u)}\over {\|u\|_X^2}}\geq {{\psi(v)}\over {\|v\|_X^2}}\geq
{{\eta}\over {\lambda_1}}\ .\eqno{(15)}$$
Now, $(13)$ is obtained from $(15)$ passing to the limit for $\eta$ tending to $\sigma$. Now, fix
$\lambda\in \left ]{{\lambda_1}\over {2\sigma}},{{\lambda_1}\over {2\max\{\rho,0\}}}\right [$. Then,
from $(10)$ and $(13)$, we obtain
$$\limsup_{\|u\|_X\to +\infty}{{\psi(u)}\over {I(u)}}<{{1}\over {\lambda}}<\sup_{u\in X\setminus \{0\}}{{\psi(u)}\over {I(u)}}\ .$$
This readily implies that $\theta^*<\lambda<\tilde\theta$ and the conclusion is directly provided by Theorem 1.\hfill $\bigtriangleup$\par
\medskip
COROLLARY 1. - {\it Let the assumptions of Theorem 2 be satisfied and let $\lambda\in \left ]{{\lambda_1}\over {2\sigma}},{{\lambda_1}\over {2\max\{\rho,0\}}}\right [$
satisfy
$$\sup_{\xi\in {\bf R}}(\lambda g(\xi)-F(\xi)f(\xi))\xi\leq 0\ .\eqno{(16)}$$
Then, for every convex set $S\subseteq L^{\infty}(\Omega)$ dense in $L^2(\Omega)$, there exists
$\alpha\in S$ such that the problem
$$\cases{-\Delta u=\alpha(x)f(u)+\lambda g(u) & in $\Omega$ \cr & \cr u=0 & on $\partial\Omega$\cr}$$
has at least three weak solutions, two of which are global minima in $H^1_0(\Omega)$ of the functional
$$u\to {{1}\over {2}}\int_{\Omega}|\nabla u(x)|^2dx-\int_{\Omega}\alpha(x)F(u(x))dx-\lambda\int_{\Omega}G(u(x))dx\ .$$}\par
\smallskip
PROOF. It suffices to observe that, in view of $(16)$, $0$ is the only weak solution of $(8)$ and then to apply Theorem 2.\hfill $\bigtriangleup$\par
\medskip
Finally, notice the following remarkable corollary of Corollary 1:\par
\medskip
COROLLARY 2. - {\it Let $q>1$, with $q<{{n+2}\over {n-2}}$ if $n\geq 3$. Let $h:{\bf R}\to {\bf R}$ be a non-negative
continuous function, with $\inf_{[0,1]}h>0$, satisfying conditions $(a)$ and $(b)$ of Theorem 2 for $f=h$.\par
Then, for every $\lambda>\lambda_1$ and for every convex set $S\subseteq L^{\infty}(\Omega)$ dense in $L^2(\Omega)$, there exists
$\alpha\in S$ such that the problem
$$\cases{-\Delta u=\alpha(x)h(u)+\lambda(u^+-(u^+)^q) & in $\Omega$ \cr & \cr u=0 & on $\partial\Omega$\cr}$$
has at least three weak solutions, two of which are global minima in $H^1_0(\Omega)$ of the functional
$$u\to {{1}\over {2}}\int_{\Omega}|\nabla u(x)|^2dx-\int_{\Omega}\alpha(x)H(u(x))dx-\lambda\int_{\Omega}\left ({{1}\over {2}}|u^+(x)|^2-{{1}\over {q+1}}|u^+(x)|^{q+1}\right )dx\ .$$}\par
\smallskip
PROOF. Fix $\lambda>\lambda_1$. Notice that, since $\inf_{[0,1]}h>0$, the number
$$\gamma:=\inf_{\xi\in ]0,1]}{{H(\xi)h(\xi)}\over {\xi}}$$
is positive. Now, we are going to apply Corollary 1 taking
$$f(\xi)=\sqrt{{{\lambda}\over {\gamma}}}h(\xi)$$
and 
$$g(\xi)=\xi^+-(\xi^+)^q\ .$$
 Of course (with the notations of Theorem 2), 
$\rho=0$ and $\sigma={{1}\over {2}}$.  Since $f$ in non-negative, $Ff$ is so in $[0,+\infty[$ and non-positive in $]-\infty,0]$. Therefore, $(16)$ is
satisfied for all $\xi\in {\bf R}\setminus [0,1]$ since $g$ has the opposite sign of $Ff$ in that set. Now, let $\xi\in ]0,1]$. We have
$${{F(\xi)f(\xi)}\over {\xi}}={{\lambda}\over {\gamma}}{{H(\xi)h(\xi)}\over {\xi}}\geq \lambda\geq \lambda(1-\xi^{q-1})$$
which gives $(16)$. Now, let $S$ be any convex set $S\subseteq L^{\infty}(\Omega)$ dense in $L^2(\Omega)$. Then, the set
$\sqrt{{{\gamma}\over {\lambda}}}S$ is convex and dense in $L^2(\Omega)$ and the conclusion follows applying Corollary 1 with this set.\hfill $\bigtriangleup$

\bigskip
\bigskip
{\bf Acknowledgement.} The author has been supported by the Gruppo Nazionale per l'Analisi Matematica, la Probabilit\`a e 
le loro Applicazioni (GNAMPA) of the Istituto Nazionale di Alta Matematica (INdAM) and by the Universit\`a degli Studi di Catania, ``Piano della Ricerca 2016/2018 Linea di intervento 2". \par
\vfill\eject
\centerline {\bf References}\par
\bigskip
\bigskip
\noindent
[1]\hskip 5pt B. RICCERI, {\it Multiplicity of global minima for
parametrized functions}, Rend. Lincei Mat. Appl., {\bf 21} (2010),
47-57.\par
\smallskip
\noindent
[2]\hskip 5pt B. RICCERI, {\it A strict minimax inequality criterion and some of its consequences}, Positivity, {\bf 16} (2012), 455-470.\par
\smallskip
\noindent
[3]\hskip 5pt B. RICCERI, {\it A range property related to non-expansive
operators}, Mathematika, {\bf 60} (2014), 232-236.\par
\smallskip
\noindent
[4]\hskip 5pt B. RICCERI, {\it Singular points of non-monotone potential operators}, J. Nonlinear
Convex Anal., {\bf 16} (2015), 1123-1129.\par
\smallskip
\noindent
[5]\hskip 5pt B. RICCERI, {\it On a minimax theorem: an improvement, a new proof and an overview of its applications},
Minimax Theory Appl., {\bf 2} (2017), 99-152.\par
\smallskip
\noindent
[6]\hskip 5pt B. RICCERI, {\it Miscellaneous applications of certain minimax theorems II}, Acta Math. Vietnam., {\bf 45} (2020), 515-524.\par
\smallskip
\noindent
[7]\hskip 5pt B. RICCERI, {\it An alternative theorem for gradient systems}, preprint.\par
\smallskip
\noindent
[8]\hskip 5pt E. ZEIDLER, {\it Nonlinear functional analysis and its
applications}, vol. III, Springer-Verlag, 1985.\par
\bigskip
\bigskip
\bigskip
\bigskip
Department of Mathematics and Informatics\par
University of Catania\par
Viale A. Doria 6\par
95125 Catania, Italy\par
{\it e-mail address}: ricceri@dmi.unict.it

\bye